\documentstyle[art12,amsfonts]{article}
\newcommand{\ns}{\normalshape}\newcommand{\nvs}{\vspace{-3pt}}
\parskip3pt\newcommand{\n}{\noindent}
\mathsurround1pt\newcommand{\k}{\kern}\newcommand{\vs}{\vspace{.3in}}
\newcommand{\ds}{\displaystyle}\newcommand{\ts}{\textstyle}
\newcommand{\sm}{\medbreak}
\renewcommand{\a}{&\k-3pt=\k-3pt&\ds}\renewcommand{\sp}[1]{{}\!^{#1}}
\newcommand{\chii}{\raise2pt\hbox{$\chi$}}\renewcommand{\P}{P}
\newcommand{\kap}{\hbox{\large$\kappa$}}
\newcommand{\la}{\lambda}
\newcommand{\Ga}{\Gamma}\newcommand{\Si}{\Sigma}
\newcommand{\si}{\sigma}\newcommand{\ga}{\gamma}
\newcommand{\R}{{\Bbb R}}\newcommand{\Z}{{\Bbb Z}}\newcommand{\C}{{\Bbb C}}
\newcommand{\CP}{\Bbb{CP}}\renewcommand{\O}{{\cal O}}
\newcommand{\F}{{\cal F}}\newcommand{\G}{{\cal G}}\newcommand{\M}{{\cal M}}
\newcommand{\qk}{quaternion-\ka}\newcommand{\ka}{K\"ahler }
\newcommand{\ch}{\hbox{\ns ch}}\def\qed{\phantom{i}\hfill QED\bigbreak}

\renewcommand{\u}{u}\renewcommand{\v}{v}\newcommand{\el}{\ell}
\newcommand{\op}{\oplus}\newcommand{\ot}{\otimes}
\newcommand{\lra}{\longrightarrow}
\newcommand{\ext}[1]{\raise1pt\hbox{$\ts\bigwedge$}\k-1pt^{#1}\k-1pt}
\newcommand{\ect}[1]{\raise1pt\hbox{$\ts\bigwedge$}_0^{#1}}
\newcommand{\sym}[1]{\raise1pt\hbox{$\textstyle\bigodot\nolimits$}
\k-.8pt^{#1}\k-.5pt} \newcommand{\Gr}{{\Bbb G}{\ns r}}
\newcommand{\scm}[1]{\raise1pt\hbox{$\textstyle\bigodot\nolimits$}_0^{#1}}
\newcommand{\ba}{\begin{array}}\renewcommand{\bar}{\begin{array}{rcl}}
\newcommand{\rf}[1]{(\ref{#1})}\newcommand{\be}{\begin{equation}}
\newcommand{\ea}{\end{array}}\newcommand{\ee}[1]{\label{#1}\end{equation}}

\newcommand{\frs}[2]{\hbox{\large$\ts\frac{#1}{#2}$\normalsize}}
\renewcommand{\ss}[1]{\vspace{30pt}\n\large{\bf#1}\setcounter{enumi}{-1}
\setcounter{equation}0\addtocounter{enumi}1\normalsize\vspace{15pt}}
\outer\def\pro#1#2\par{\bigbreak\noindent{\bf#1\enspace}{\sl#2}\par\bigbreak}
\newcommand{\pf}{\n{\sl Proof. }}\newcommand{\re}{\n{\sl Remark. }}
\begin{document}

\centerline{\Large\bf The Twistor Transform of a Verlinde formula}
\normalsize\vs\centerline{\large S.$\,$M.\ Salamon}\vs

\ss{Introduction}

Let $\Si$ be a compact Riemann surface of genus $g$. The moduli space
$\M_g=\M_g(2,1)$ of stable rank 2 holomorphic bundles over $\Si$ with
fixed determinant bundle of degree 1 is a smooth complex
$(3g-3)$-dimensional manifold \cite{Ses}. The anticanonical bundle of
$M$ is the square of a holomorphic line bundle $L$, some power of
which embeds $\M_g$ into a projective space. The dimensions of the
vector spaces $H^0(\M_g,\O(L^{m-1}))$ of holomorphic sections of
powers of $L$ are known to be independent of the choice of complex
structure on $\Si$, and are given by the formula \be
h^0(\M_g,\O(L^{m-1}))=-m^{g-1}\sum_{i=1}^{2m-1}(-1)^i\hbox{cosec}^{2g-2}
(\frac {i\pi}{2m})\ee{V} predicted by Verlinde \cite{V}. This is
closely related to the structure of the cohomology ring of $\M_g$, and
a number of independent proofs and generalizations of \rf{V} are now
known. Below we shall follow closely the approach of Szenes \cite{Sz}.

In the case in which $\Si$ is a hyperelliptic surface, and is
therefore a 2-fold branched covering of $\CP^1$, Desale and Ramanan
\cite{DR} exhibit $\M_g$ as a complex submanifold of the flag manifold
$\F_g=SO(2g+2)/(U(g-1)\times SO(4))$. As explained in
\cite{Sz} this reduces verification of \rf{V} to certain
$SO(2g+2)$-equivariant calculations. Our contribution is to observe
that $\F_g$ is the twistor space of the real oriented Grassmannian
$\G_g= SO(2g+2)/(SO(2g-2)\!\times\!SO(4))$ in the sense of
\cite{Br,BR} for all $g\ge3$. This enables us to relate the cohomology
of the symmetric space $\G_g$ directly to the cohomology of $\M_g$,
and we obtain a set of generators for the latter which may be compared
to the universal ones described in \cite{N2,Th,D}. As a feasibility
study, we illustrate the theory in the present paper for the case
$g=3$ which is worthy of special attention since the fibration
$\F_3\to\G_3$ encapsulates the quaternionic structure of the base
space in a manner first identified by Wolf \cite{Wo}.

In the first section we investigate the cohomology of $\G_3=SO(8)/
(SO(4)\!\times\!SO(4))$. Using its quaternionic spin structure, we
prove that the odd Pontrjagin classes of $\G_3$ vanish, and that its
$\hat A$ class simplifies remarkably. In the second section we recover
Ramanan's description \cite{R} of the Chern ring of $\M_3$ in the
context of the natural mapping $\M_3\to\G_3$, enabling
$h^0(\M_3,L^{m-1})$ to be computed rapidly. Whilst this provides only
a particularly simple instance of
\rf{V}, results of the third section identify $H^0(\M_3,L^k)$ with a
virtual representation of $SO(8)$ that also arises from the kernels of
coupled Dirac operators on $\G_3$. Similar techniques can in theory be
applied to higher genus cases, and formulae such as $p_1\sp g=0$ on
$\M_g$ \cite{Th,Ki} may be expected to interact with properties of
$\G_g$ such as the constancy of the elliptic genera considered in
\cite{W,HS}.\smallbreak

The material below was presented at the conference `Differential
Geometry and Complex Analysis' in Parma in May 1994, and the author is
grateful to the organizers of that event.

\ss{1. Grassmannian cohomology}

{}From now on we denote by $\G$ the Grassmannian
\be \G_3=\frac{SO(8)}{SO(4)\times SO(4)} \ee{iso}
that parametrizes real oriented 4-dimensional subspaces of $\R^8$. Let
$W$ denote the tautological real rank 4 vector bundle over $\G$, and
$W^\perp$ its orthogonal complement in the trivial bundle over $\G$
with fibre $\R^8$. The bundles $W$ and $W^\perp$ arise from the
standard representations of the two $SO(4)$ factors constituting the
isotropy subgroup in \rf{iso}, and it follows that
\be    T\G\cong W\ot W^\perp.\ee{tp}  The $SO(8)$-invariant
Riemannian metric on $\G$ determines an isomorphism $W\cong W^*$ of
vector bundles.

The decomposition \rf{tp} may be refined by lifting the $SO(4)$
structure of $W$ to $\hbox{\it Spin}(4)\cong SU(2)\times SU(2)$ on a
suitable open dense subset $\G'$ of $\G$. This procedure is one that
is familiar from the study of Riemannian 4-manifolds, and
\[     W_\C \cong U\ot_\C V,\]
where $U$ and $V$ are each complex rank 2 vector bundles over $\G'$.
The resulting isomorphism
\be  (T\G)_\C\cong U\ot(V\ot W^\perp_\C)  \ee{quat}
reflects the fact that $\G$ is a \qk manifold \cite{Wo,S}.
In \rf{quat}, $U$ may be thought of as a quaternionic line bundle
(usually called $H$), and its cofactor $V\ot W^\perp_\C$ (usually
called $E$) has structure group $SU(2)\times SO(4)$ extending to
$Sp(4)$.

The Betti numbers of a \qk $4n$-manifold of positive scalar curvature
satisfy $b_{2k+1}=0$ for all $k$ and $b_{2k-4}\le b_{2k}$ for $k\le
n+1$. They are also subject to the linear constraint of \cite{LS}
which for $n=4$ takes the form
\[  3(b_2+b_4)=1+b_6+2b_8.\]
This is well illustrated by $\G$, which has
Poincar\'e polynomial
\[   P_t(\G)=1+3t^4+4t^8+3t^{12}+t^{16},\]
and is the only real Grassmannian to have $b_4>2$. (These facts may be
deduced from \cite[chapter XI]{GHV}.) We shall in fact only be
concerned with the subring generated by the Euler class $e=e(W)$ and
the first Pontrjagin class $f=p_1(W)$.

Although the classes $e$ and $f$ are very natural, it will ultimately
be more convenient to consider
\[ u=-c_2(U),\quad v=-c_2(V).\]
Because of the $\Z_2$-ambiguity in the definition of $U,V$, the
classes $u,v$ are not integral, but the symmetric products $\sym2U,
\sym2V$ are globally defined so $4u,4v$ belong to $H^4(\G,\Z)$. If we
write formally $4u=\ell^2$ then
\be \ch(U)=e^{\ell/2}+e^{-\ell/2}=2+u+\frs1{12}u^2+\frs1{360}u^3+
\frs1{20160}u^4.\ee{ch}
The class $\ell$ is given geometrical significance by the splitting
\rf{formal}. An analogous expression to \rf{ch} holds for $\ch(V)$, and
from $\ch(W_\C)= \ch(U)\ch(V)$, we obtain \be\ba{c} e = u-v,\\ f =
2(u+v).\ea\ee{comp} We may add that $p_2(W)=c_4(W_\C)=(u-v)^2$
confirming the well-known relation
\be  p_2(W)=e^2. \ee{p2} Moreover, the space $H^4(\G,\Z)$ is generated by
$e,f$ together with $e(W^\perp)$ \cite{MS}.

\pro{1.1 Proposition} Evaluation on the fundamental cycle $[\G]$ yields
\nvs\[\ba c e^4=2=e^2f^2,\quad e^3f=0=ef^3,\quad f^4=4;\\[5pt] u^4=
\frs{21}{64}=v^4,\quad u^3v=-\frs7{64}= uv^3,\quad u^2v^2=\frs5{64}.\ea\]

We shall deduce these Schubert-type relations from a description of
the total Pontrjagin class and the $\hat A$ class
\[\bar  \P(T\G)\a1+\P_1+\P_2+\P_3+\P_4,\\[3pt] \hat A(T\G)\a1+\hat A_1+
\hat A_2+\hat A_3+\hat A_4\ea\] of the tangent bundle \rf{tp} of $\G$.
(Upper case $\P_i$'s are used to prevent a future clash of notation.)
The classes $\hat A_i$, $1\le i\le 4$ are determined in terms of the
$\P_i$ in the usual way \cite{H}, and

\pro{1.2 Proposition} $\P_1=0=\P_3$ and $\hat A(\G)=1-\frs1{240}f^2$.

\n{\sl Proof of both propositions.} It is easy to check that, in the
presence of \rf{comp}, the two sets of equations of Proposition~1.1
are equivalent.  The equalities $u^4=v^4$ and $u^3v=uv^3$ are
immediate from the symmetry between $U$ and $V$, and these are
equivalent to $ef^3=0=ef^3$. Using \rf{p2}, we have
\be  \ch(W_\C)= 4+f+\frs1{12}(-2e^2+f^2)+\frs1{360}(-3e^2f+f^3)+
\frs1{20160}(2e^4-4e^2f^2+f^4).\ee{chW}
{}From \rf{tp} and \rf{chW},
\be\bar \ch(T\G)_\C\a(\ch\,W_\C)(8-\ch W_\C)\\[5pt]\a
    16-f^2+\frs16(2e^2f-f^3)+\frs1{720}(-20e^4+32e^2f^2-9f^4).
\ea\ee{chT} In particular $\P_1=0$, and so we also have
\be \ch(T\G)_\C=16-\frs16\P_2+\frs1{120}\P_3+\frs1{10080}(\P_2\sp2-2\P_4).
\ee{chTT} Comparing \rf{chT} and \rf{chTT} gives
\be \P_2=6f^2,\quad\P_3=20(2e^2f-f^3),\quad\P_4=140e^4-224e^2f^2+81f^4.
\ee{PPP}\sm

The remainder of the proof is based on the following less obvious
facts.\sm\par\n(i) $\G$ is a spin manifold (see forward to \rf{+-})
carrying a metric of positive scalar curvature. Therefore its $\hat A$
genus\be \hat A_4=\frs1{2^{16}3^45^27}(762\P_1\sp4 -
1808\P_1\sp2\P_2 + 416\P_2\sp2 + 1024\P_1\P_3 - 384\P_4)\ee{A4}
vanishes. Thus \be 0 = 416(6f)^2 - 384(140e^4-224e^2f^2+81f^4) =
5376(-10e^4+16e^2f^2 -3f^4).\ee1\sm

\n(ii) The dimension $d$ of the isometry group of any \qk 16-manifold
with positive scalar curvature is given by
\[  d = 7-\frs83\P_1u^3+64u^4\]
\cite[page~170]{S}. In the present case, $d=\dim SO(8)=28$ and we obtain
\be  21 = 64u^4 = \frs14(16e^4+24e^2f^2+f^4).\ee2\sm

\n(iii) On any compact \qk $4n$-manifold $M$ with positive scalar
curvature and $n>2$, the index
\[ \hat A(M,\sym2U)=\left<\ch(\sym2U)\hat A,[M]\right>,\]
vanishes; this is a consequence of \cite[Corollary~6.7]{S}
which is explained in \cite{LS}. Given that
\[\ba c\ch(\sym2U) = 3+4u+\frs43u^2+\frs8{45}u^3+\frs4{315}u^4,\\[8pt]
\hat A = 1-\frs1{24}\P_1-\frs1{2^5 3^2 5}\P_2-\frs1{2^6 3^3 5^1 7}\P_3=
1-\frs1{240}f^2+\frs1{1008}(2e^2f-f^3),\ea\]
and $u=(2e+f)/4$, it follows that
\be  24e^4 - 26e^2f^2 + f^4 = 0.\ee3\sm
Proposition~1.1 now follows from \rf1,\rf2,\rf3, and it only
remains to prove that $\P_3=0$. Because of the symmetry between $W$
and $W^\perp$, it suffices to prove that $\P_3e=0=\P_3f$, but this
follows from \rf{PPP} and Proposition~1.1.\qed

\re The vanishing of $\hat A_4$ and \rf1 above is in fact equivalent
to the vanishing of the index $\hat A(M,T)$ of the Dirac operator
coupled to the tangent bundle (see \rf{AS}), essentially the so-called
Rarita-Schwinger operator. This index is known to be equivariantly
constant on any spin manifold with $S^1$ action \cite{W}, and always
vanishes in the homogeneous setting \cite{HS}.

\ss{2. The flag manifold and moduli space}

We denote by $\F$ the complex 9-dimensional homogeneous space
\be   \F_3=\frac{SO(8)}{U(2)\times SO(4)}  \ee{coset}
that parametrizes complex 2-dimensional subspaces $\Pi$ of $\C^8$ that
are isotropic with respect to a standard $SO(8)$-invariant bilinear
form. It has a complex contact structure that was studied in \cite{Wo}
and exhibits it as the twistor space of $\G$ in the sense of
\cite{S}. Projecting $\Pi$ to a real 4-dimensional subspace of $\R^8$
determines an $SO(8)$-equivariant mapping $\pi\colon\F\to \G$, and
each fibre of $\pi$ is isomorphic to $SO(4)/U(2)$ and defines a
rational curve in the complex manifold $\F$.

{}From standard facts about twistor spaces \cite{Be,S,PS}, one knows
that $\hbox{Pic}(\F)$ is generated by a holomorphic line bundle $L$ on
$\F$ such that\sm\par\n(i) the restriction of $L$ to each fibre
$\pi^{-1}(x)\cong\CP^1$ equals $\O(2)$;\par
\n(ii) $L^5$ is isomorphic to the anticanonical bundle $\kap^{-1}$ of
$\F$. \sm\par\n The line bundle $L$ admits a square root over an open
set $\G'$ of $\G$ on which $U$ and $V$ are defined, there is a
$C^\infty$ isomorphism
\be  \pi^*U\cong L^{1/2}\op L^{-1/2}. \ee{formal}
Let $\ell$ denote the fundamental class $c_1(L)$ in $H^2(\F,\Z)$. From
the Leray-Hirsch theorem, there is an identity
$(\ell/2)^2+\pi^*c_2(U)=0$ of real cohomology classes. In terms of
integral classes, and omitting $\pi^*$, \be \ell^2=4u.\ee{LH}\sm

In the notation of the Introduction, let $\M=\M_3$. Szenes exhibits
the latter as the zero set of a non-degenerate holomorphic section
$s\in H^0(\F,\O(\si^*))$, where $\si=\sym2\tau$ and $\tau$ denotes the
tautological rank 2 complex vector bundle acquired from the embedding
$\F\subset\Gr_2(\C^8)$. (Such a section $s$ corresponds to a quadratic
form on $\C^8$, but we shall not mention this again until the end of
Section~3.)  From the coset description \rf{coset}, it follows that
\be \tau\cong L^{-1/2}\ot\pi^*V;\ee{taut} the right-hand side is well
defined on $\F$, even though the individual factors only make sense
locally (for example on $\pi^{-1}(\G')$). Since $V\cong V^*$, we have
$\si^*\cong L\ot\pi^*\sym2V$. The resulting holomorphic structure on
$\pi^*\sym2V$ coincides with that induced in a standard way from the
fact that $\sym2V$ has a self-dual connection on the \qk manifold
$\G$, in the sense of \cite{MCS}. In particular, $\pi^*\sym2V$ is
trivial over each fibre $\pi^{-1}\cong\CP^1$. From now on we shall
write $\sym2V$ in place of $\pi^*\sym2V$, and often omit tensor
product signs.

The cohomology classes $\ell,u,v$ may be pulled back from both $\G$
and $\F$ to $\M$, and we shall denote the resulting elements of
$H^i(\M,\R)$ by the same symbols.

\pro{2.1 Proposition} On $\M$, $3\u^2+10\u\v +3\v^2=0$, and
evaluation on $[\M]$ yields\nvs\[ \u^3=\frs72=-\v^3,\quad\u\v^2=
\frs32=-\u^2\v.\]

\pf The submanifold $\M$ of $\F$ is Poincar\'e dual to the Euler class
$c_3(\si^*)$, which is readily computed from the formula
$\ch(\si^*)=e^\ell\ch(\sym2V)$ (see \rf{taut}) and equals
$4\ell(u-v)$. Then, for example,
\[ \left<\u^3,[\M]\right>=\left<u^3\,c_3(\si^*),[\F]\right>=
 \left<4\ell(u^4-u^3v),[\F]\right>=8\left<u^4-u^3v,[\G]\right>=\frs72,\]
the last equality from Proposition~1.1. The evaluation of $\u^2\v$,
$\u\v^2$ and $\v^3$ follows in exactly the same way.

Since $H^8(\M,\R)\cong H^4(\M,\R)$ is 2-dimensional \cite{N}, there
must be a non-trivial linear relation $a\u^2+b\u\v+c\v^2=0$. The
solution $(a=c)/b=3/10$ can be found by multiplying the left-hand side
by $u$ and $v$ in turn.\qed

The next result gives an independent derivation of the characteristic
ring in the context of the twistor fibration $\F\to\G$.

\pro{2.2 Proposition} The Chern and Pontrjagin classes of $\M$ are given
by \nvs\[\ba c c_1=2\el,\quad c_2=4(3\u+\v),\quad c_3=8\el\u,\quad
c_4=-\frs{112}3\u\v,\quad c_5=c_6=0;\\[4pt] p_1=-8(\u+\v),\quad p_2=
\frs38p_1\sp2,\quad p_3=0.\ea\]

\pf It is known \cite{S} that the fibration $\pi$ gives a $C^\infty$
splitting of the holomorphic tangent bundle of $\F$:
\[ T^{1,0}\F\cong L\,\op\,L^{1/2}(V\ot W^\perp_\C).\]
Combining this with the isomorphism
\[ T^{1,0}\F|_{\M}\cong T^{1,0}\M\op(L\,\sym2V)|_{\M},\] we obtain
\[\bar \ch(T^{1,0}\M)\a e^{\el}+e^{\el/2}\ch(V\,W_\C^\perp)-e^\el
\ch(\sym2V) \\[3pt]\a e^{\el}(1+e^{-\el/2}\ch V(8-\ch W_\C)-
\ch(\sym2V)).\ea\] This yields the required expressions for
$c_1,c_2,c_3$. We also get $c_4=28(\u+\v)^2$ which reduces to
$-112\u\v/3$ from Proposition~2.1. We next obtain $c_5=-32\el v(u+v)$,
so that $c_5\ell=0$ and the vanishing of $c_5$ follows from the fact
that $H^2(\M,\R)$ is 1-dimensional \cite{N}. Finally, all these
equalities combine to yield \[
c_6=\frs13(504u^3+2824u^2v+1928uv^2+120v^3),\] and Proposition~2.1
implies that $c_6=0$. The Pontrjagin classes $p_i$ of $\M$ are now
determined from the Chern classes by the usual relations.\qed

\re The cohomology ring and Chern classes of $\M$ were
computed in \cite[Theorem~4]{R}, and comparison with that shows that
\[h=\el,\quad \nu=\frs12(3u+v).\] In general, it is known that the total
Pontrjagin class of $\M_g$ equals $(1+\frs1{2g-2}p_1)^{2g-2}$
\cite{N2}. Moreover, $p_1\sp g=0$ \cite{Ki,Th} and $c_i=0$ if $i>2g-2$
\cite{Gie}.\sm

The above enable the dimension $d_k$ of $H^0(\M,\O(L^k))$ to be
computed quickly. For this purpose it is convenient to set $k=m-1$.

\pro{2.3 Theorem} $d_{m-1}=\frs1{45}m^2(11+20m^2+14m^4)$.

\pf Given that $c_1(T^{1,0}F)=2\ell$, the Todd class
$\hbox{td}(T^{1,0}\M)$ of $\M$ equals \[ e^\el\hat A(T\M)=
e^\el\left[1-\frs1{24}p_1+\frs1{2^73^25}(7p_1\sp2-4p_2)\right].\]
Using Propositions~2.1, 2.2 and the Riemann-Roch theorem, we obtain
\[\bar d_{m-1}\a\left<e^{m\el}(1+\frs13(u+v)-\frs{11}{135}\u\v),
[\M]\right>\\[8pt]\a -\frs{22}{135}m^2\u^2\v+\frs29m^4(\u^3+\u^2\v)
+\frs4{45}m^6\u^3,\ea\] and the result follows.\qed

\ss{3. Equivariant indexes}

In this section, we begin by considering the Dirac operator over the
Grassmannian $\G$. Recall from \rf{quat} that the quaternionic
structure of $\G$ is characterized by the vector bundles $H=U$ and
$E\cong V\>W_\C$ (juxtaposition denotes tensor product). For $p\ge4$,
the exterior power $\ext pE$ contains a proper subbundle $\ect pE$
with the property that $\ext pE\cong\ect pE\op \ext{p-2}E$ and, as
described in \cite{BS}, the total spin bundle $\Delta$ of $\G$
decomposes as $\Delta_+\op\Delta_-$ where
\be\bar \Delta_+ &\cong& \sym4U\ \op\ \sym2U\,\ect2E\ \op\ \ect4E,\\
\Delta_- &\cong& \sym3U\,E\ \op\ U\,\ect3E.\ea\ee{+-}
The fact that all the summands on the right-hand side are globally
defined confirms that $\G$ is spin, though we shall not in fact need
the decompositions \rf{+-}.

Now let $X$ be any other complex vector bundle over $\G$. The choice
of a connection on $X$ allows one to extend the Dirac operator on $\G$
to an elliptic operator
\[ D_X\colon\Ga(\Delta_+\>X)\lra\Ga(\Delta_-\>X).\]
The index of this coupled Dirac operator is by definition $\dim(\ker
D_X)-\dim(\hbox{coker}\,D_X)$. This extends to a homomorphism
$K(\G)\to\Z$, so that the index of $D_X$ is also defined when $X$ is a
virtual vector bundle. The Atiyah-Singer index theorem \cite{AS}
asserts that the index of $D_X$ equals
\be \hat A(\G,X)=\left<\ch(X)\hat A(T\G),[\G]\right>\ee{AS}
In our situation, this fact is closely related to the Riemann-Roch
theorem on $\F$ which provides the following interpretation of $d_k$.

\pro{Theorem 3.1} Let $X_k=\sym{2k+4}U-\sym{2k+2}U\sym2V+\sym{2k}U\sym2V-
\sym{2k-2}U$, $k\ge1$. Then $d_k=\hat A(\G,X_k)$.

\pf Let $\si$ denote the rank 3 vector bundle $\sym2\tau$ as above, and
let $(k)$ denote the operation of tensoring with $L^k$. The
description of $\M$ as the zero set of a section of
$\si^*\cong\sym2V(1)$ provides a Koszul complex
\[ 0\to\O_{\F}(\ext3\si(k))\to\O_{\F}(\ext2\si(k))\to\O_{\F}(\si(k))\to
\O_{\F}(k)\to\O_{\M}(k)\to0,\]
or equivalently,
\[0\to\O_{\F}(k-3)\to\O_{\F}(\sym2V(k-2))\to\O_{\F}(\sym2V(k-1))\to
\O_{\F}(k)\to\O_{\M}(k)\to0.\]
It follows that \be \chii(\M,\O(k))= a_k-b_{k-1}+b_{k-2}-a_{k-3},\ee4
where \be a_k=\chii(\F,\O(k)),\quad
b_k=\chii(\F,\O(\sym2V(k))).\ee{ab} These holomorphic Euler
characteristics may be computed using the Riemann-Roch theorem and the
cohomological version
\cite[7.2]{S} of the twistor transform; the result is
\be a_k=\hat A(\G,\sym{2k+4}U),\quad b_k=\hat A(\G,\sym{2k+4}U\sym2V).
\ee{ab2} Finally, Proposition~2.2 implies that the canonical bundle
$\kap(\M)$ is isomorphic to $L^{-2}$, so by Serre duality and Kodaira
vanishing, $H^i(\M,\O(k))=0$ for all $i\ge1$ and $k\ge-1$. In
particular, $\chii(\M,\O(k))=\dim H^0(\M,\O(k))$ for all $k\ge-1$, and
the theorem now follows from \rf{4}.\qed

The isometry group $SO(8)$ of $\G$ acts naturally on the cohomology
groups over $\F$ of the sheaves $\O(k),\,\O(\sym2V(k))$ considered
above. The integers $a_k,\,b_k$ and \[d_k=a_k-b_{k-1}+b_{k-2}-a_{k-3}
\] are therefore the dimensions of certain virtual $SO(8)$-modules,
and we identify these shortly.

Let $V(\ga)$ denote the complex irreducible representation of $SO(8)$
with dominant weight $\ga$, where $\ga=(\la_1,\la_2,\la_3,\la_4)$ with
$\la_1\ge\la_2\ge \la_3\ge\la_4\ge0$. We adopt standard coordinates so
that $V(1,0,0,0)=\C^8$ is the fundamental representation, and
$V(1,1,0,0)=\frak{so}(8,\C)$ is the complexified adjoint
representation.

\pro{3.2 Proposition} Let $A_k=V(k,k,0,0)$ and $B_k=V(k+1,k-1,0,0)$.
Then $a_k=\dim A_k$ and $b_k=\dim B_k$.

\pf The Weyl dimension formula states that
\[   \dim(V(\ga))=\prod_{\alpha\in R_+}\frac{\left<\alpha,d+\ga\right>}
{\left<\alpha,d\right>},\] where $R_+$ denotes the set of positive
roots and $d$ is half of their sum. With the above coordinates,
\newcommand{\y}{\hbox{$-$}\kern-.5pt} \[\ba{l}
R_+=\{(1,1,0,0),(1,0,1,0),(1,0,0,1),(0,1,1,0),(0,1,0,1),(0,0,1,1),
\\\hspace{40pt} (1,\y1,0,0),(1,0,\y1,0),(1,0,0,\y1),(0,1,\y1,0),(0,1,0,\y1),
(0,0,1,\y1)\}, \ea\] $d=(3,2,1,0)$ and we obtain
\[\bar \dim A_k\a\frs1{4320}(k+1)(k+2)^3(2k+5)(k+3)^3(k+4),
\\[6pt]\dim B_k\a\frs1{1440}k(k+1)^2(k+2)(2k+5)(k+3)(k+4)^2(k+5).\ea\]

We claim that the right-hand sides are equal to $a_k$ and $b_k$
respectively. It follows from \rf{ab} that $a_k$ and $b_k$ are
polynomials in $k$ of degree 9, and by Serre duality,
\be  a_{-k}=-a_{k-5},\quad b_{-k}=-b_{k-5},\qquad k\in\Z.\ee{Serre}
By \rf{Serre} and suitable vanishing theorems \cite{BE}, $a_k=0=b_k$
for $k=-4,-3,-\frs52,-2,-1$. In addition, $\F$ has Todd genus
$a_0=1=-a_{-5}$, and $b_0=0=b_{-5}$. Accordingly,
\[\bar   a_k\a\frs1{4320}(k+1)(k+2)(2k+5)(k+3)(k+4)\tilde a_k,\\[6pt]
 b_k\a\frs1{1440}k(k+1)(k+2)(2k+5)(k+3)(k+4)(k+5)\tilde b_k.\ea\]
where $\tilde a_k$ is a quartic polynomial in $k$ with $\tilde
a_0=36$ and $\tilde b_k$ is quadratic in $k$.

Let $n=2k+4$. The formulae \rf{ab2} involve $\ch(\sym nU)=f(n)$, where
\be f(x)=
\frac{e^{(x+1)\ell/2}-e^{-(x+1)\ell/2}}{e^{\ell/2}-e^{-\ell/2}}\ee{quo}
(see \rf{formal}). To evade an explicit calculation of $\ch(\sym nU)$,
we exploit the following formulae which are easily deduced from
\rf{quo}.

\pro{3.3 Lemma} $\ds f'(0)=\frac{\ell/2}{\hbox{\ns tanh}(\ell/2)},\quad
f''(0)=u$.

\n The right-hand side of the first equation is the
series used in the definition of Hirzebruch's L-genus, and using
\rf{LH} can be rewritten as
\[ \bar \ds \frac d{dn}\Big|_{n=0}\ch(\sym nU) \a 1 -
 \sum_{j\ge1}(-1)^j\frac{2^{2j}B_j}{(2j)!} u^{2j}\\[8pt]\a
\frs12(1-\frs13u-\frs1{45}u^2+\frs2{945}u^3-\frs1{4725}u^4),\ea\]
where $B_j$ are the Bernoulli numbers \cite{H}. From above, we obtain
\[ \ds\frac d{dk}\Big|_{k=-2}a_k= \frs1{270}(u^4+2u^3v+u^2v^2)-
\frs1{4725}u^4=0=\frac{d^2}{dk^2}\Big|_{k=-2}a_k.\]
It follows that $\tilde a_k$ is divisible by $(k+2)^2$, and by Serre
duality by $(k+3)^2$. We obtain $\tilde a_k=(k+2)^2(k+3)^2$.  The
identification $\tilde b_k=(k+1)(k+4)$ is similar, and proceeds using
a less-enlightening version of Lemma~3.3; we omit the details.\qed

The following table displays some of the above dimension functions in
terms of $k$.

\small\[ \ba{|c||c|c|c|c|c|c|c|c|c|}\hline k&0&1&2&3&4&5&6&7&8\\\hline a_k
&\,1\,&28&300&1925&8918&32928&102816&282150&698775\\\hline b_k
&0&35&567&4312&21840&85050&274890&772464&1945944\\\hline
d_k &1&28&265&1392&5145&15100&37681&83392&168273\\\hline
\ea\]\normalsize\bigbreak\vspace{.2in}

Applying Serre duality and Kodaira vanishing over $\F$, recalling that
$\kap(\F)\cong L^{-5}$, shows that there is in fact an
$SO(8)$-equivariant isomorphism $A_k\cong H^0(\F,\O(k))$. In
particular, $A_1$ may be identified with both the space of holomorphic
sections of $L$ and the Lie algebra $\frak{so}(8,\C)$ of infinitesimal
automorphisms of the contact structure of $\F$. There is an associated
moment mapping $\F\to{\Bbb P}(\frak{so}(8,\C)^*)\cong {\Bbb CP}^{27}$
that identifies $\F$ with the projectivization of the nilpotent orbit
of minimal dimension
\cite{Sw}. Accordingly, the $SO(8)$-equivariant linear mapping \be
\phi_k\colon\sym k(H^0(\F,\O(1)))\lra H^0(\F,\O(k))\ee{ev} is onto for all
$k\ge1$. Indeed, $A_k$ is the irreducible summand of $\sym k A_1$ of
highest weight, and it suffices to show that the restriction of
$\phi_k$ to $A_k$ is an isomorphism. Observe that $A_k$ contains a
decomposable tensor product $\xi^{\ot k}$ for some non-zero $\xi\in
A_1$ and $\phi_k(\xi^{\ot k})$, being the $k$th power of $\xi$
regarded as a section of $L$, is also non-zero. The irreducibility of
$A_k$ and Schur's lemma establishes the claim.

A similar argument can be given to establish an $SO(8)$-equivariant
isomorphism $B_k\cong H^0(\F,\O(\sym2V(k)))$, given that
$H^i(\F,\O(\sym2V(k)))$ vanishes for all $i>0$ and $k\ge0$. One
considers the mapping
\[ \psi_k\colon H^0(\F,\O(\sym2V(1)))\ot H^0(\F,\O(k-1))\lra
H^0(\F,\O(\sym2V(k))),\] in which $H^0(\F,\O(\sym2V(1)))$ is
isomorphic to the irreducible $35$-dimensional $SO(8)$-module
$\scm2\C^8$ with highest weight $(2,0,0,0)$. The irreducible summand
of highest weight in the tensor product is isomorphic to $B_k$ and the
restriction of $\psi_k$ to this is an isomorphism.

The above arguments can be streamlined by applying more sophisticated
twistor transform machinery contained, for example, in
\cite{BE}. In particular, $A_k$ and $B_k$ are known to be
isomorphic to the respective kernels of natural twistor operators
\[  \ba l  \alpha_k\colon \sym{2k}U\lra E\>\sym{2k+1}U,\\
       \beta_k\colon \sym{2k}U\>\sym2V\lra
E\>\sym{2k+1}U\>\sym2V.\ea\] Recall that $\M$ is the zero set of an
element $s$ of the space $B_1\cong\scm2\C^8$. For suitable
hyperelliptic surfaces $\Si$, the section $s$ will be a real element;
at each point of $\G$ it then defines a section of $W\op W^\perp$,
which is a trivial bundle with fibre $\R^8$ (see \rf{iso}). In these
terms the element $\tilde s\in\ker\beta_1$ determined by $s$ is
essentially the image of $s$ by the homomorphism
\[ \sym2(W\op W^\perp)_\C\to\sym2W_\C\to\sym2U\,\sym2V\cong\hbox{Hom}(
\sym2V,\sym2U).\] This may be used to describe $\M$ as a `branched cover'
of a real subvariety of $\G$.

The Horrocks instanton bundle over $\CP^5$ discussed at the end of
\cite{MCS} provides an analogous situation in which a geometric object
is defined by a non-degenerate solution of a twistor equation over a
homogeneous space. Such situations are worthy of more systematic
investigation.

\vfil\eject

\renewcommand{\thebibliography}{\list{{\bf\arabic{enumi}.\hfil}}
{\settowidth\labelwidth{18pt}\leftmargin\labelwidth\advance
\leftmargin\labelsep\usecounter{enumi}}\def\newblock{\hskip.05em}
\sloppy \sfcode`\.=1000\relax}\baselineskip0pt
\newcommand{\bi}{\vspace{-5.6pt}\bibitem}

\ss{References}

\medskip\n Mathematical Institute, 24--29 St.\ Giles', Oxford
OX1$\,$3LB

\enddocument